\documentclass[12pt,reqno,pre,longbibliography]{amsart}
\usepackage{amsfonts}
\usepackage{graphicx}
\usepackage{amscd}
\usepackage{amsmath,amssymb}
\usepackage{epsfig}
\usepackage[usenames]{color}
\usepackage{subfigure}
\usepackage[english]{babel}

\providecommand{\U}[1]{\protect\rule{.1in}{.1in}}
\providecommand{\U}[1]{\protect\rule{.1in}{.1in}}
\allowdisplaybreaks
\newtheorem{theorem}{Theorem}

\theoremstyle{plain}
\newtheorem{acknowledgement}{Acknowledgement}

\numberwithin{equation}{section}

\DeclareMathOperator{\sech}{sech}

\DeclareMathOperator{\dn}{dn}
\DeclareMathOperator{\cn}{cn}
\DeclareMathOperator{\sn}{sn}

\begin{document}
\title[Coupled Nonlinear Schr\"odinger-Type Equations ]{Bifurcation Analysis, Modulational Instability  and Solitons Solutions for a Coupled Nonlinear Schr\"odinger-Type Equations with Variable Coefficients}
\author{José M. Escorcia }
\address{Escuela de Ciencias Aplicadas e Ingenier\'ia, Universidad EAFIT, Carrera 49 No. 7 Sur-50, Medellín 050022, Colombia.}
\email{jmescorcit@eafit.edu.co}

\date{\today }


\begin{abstract}
Understanding, predicting, and controlling physical processes often relies on the analysis of the dynamics of partial differential equations (PDEs). In this context, the present study offers an in-depth investigation into the nonlinear dynamics of a novel coupled nonlinear Schr\"odinger system with variable coefficients. To begin with, the dynamics of traveling-wave solutions is analyzed through bifurcation theory, which uncovers the presence of Jacobi elliptic functions and solitonic structures. Then, a sensitivity analysis is carried out to confirm the numerical relevance and stability of these solutions. However, when subjected to external perturbations, the system exhibits a tendency toward chaotic behavior. Finally, the results point to the emergence of modulational instability under specific configurations of the dispersion and nonlinearity parameters. \\
 

\end{abstract}

\maketitle

\section{Introduction}

The coupled nonlinear Schr$\ddot{\mbox{o}}$dinger (CNLS) equations describe various physical phenomena, including water waves \cite{Shukla2006,Zhao2016,Roskes}, Bose-Einstein condensates \cite{Deconinck2004,Santos2023,Li2025,Theocharis2003,Qi2012}, and nonlinear optics \cite{Xu2024,Triki2024,Chai2017}, to name a few. Many studies on these types of models take as their starting point the well-known Manakov system, introduced in 1974 \cite{Manakov1974}, an integrable model characterized by having equal nonlinear interactions both within each component and between them \cite{ZAKHAROV1982}. To analyze CNLS equations, a variety of powerful mathematical techniques are used: the inverse scattering transform \cite{BIONDINI2015}, the Darboux transformation \cite{WRIGHT2003}, the Hirota bilinear method \cite{Radhakrishnan1997}, and the Painlevé analysis \cite{Choudhury2006}, among others. These techniques aim to derive exact or approximate solutions, such as solitons, breathers, and rogue waves.

Recognizing the importance of coupled NLS systems in theoretical and applied sciences, generalized models have been introduced \cite{BoLi,Du,Rehab,Kevrekidis,Manganaro,Arash,Jose2024}, and considerable research has focused on exploring their nonlinear dynamics. In particular, these systems have been investigated from various angles, including bifurcation analysis, the emergence of chaotic behavior, and modulation instability \cite{Raza2023,Tang2023,Nadeem2025,Patel2021,Frisquet2015}. Bifurcation analysis provides insight into how stable states can evolve into unstable ones, often manifesting through the emergence of different types of solutions. Chaos theory emphasizes the inherent complexity and unpredictability of the system. Both concepts offer valuable insight into the natural dynamics of the system, deepening our understanding of real-world phenomena, and help to predict and effectively control them. 
 
On the other hand, modulation instability (MI), also known as Benjamin-Feir instability, is a fundamental mechanism where a slowly modulated continuous wave becomes unstable, leading to the exponential amplification of perturbations and the eventual formation of an envelope soliton train \cite{Patel2021,Rapti2004}. MI in Bose-Einstein condensates has attracted significant interest from both theoretical and experimental perspectives, as studies indicate that it may play a key role in phenomena such as condensate dephasing and localization \cite{Theocharis2003,Qi2012,Everitt2017}. It is worth noting that modulation instability also arises in various other fields, such as optical fibers \cite{Hasegawa1984}, water wave dynamics \cite{Stuhlmeier2024}, and biological systems \cite{Turing1952} .

In the present work, we introduce the following CNLS equations:
\begin{eqnarray}
\label{Eq1a}
i\psi _{t} &=&-a(t) \psi_{xx} + b(t)x^2\psi-ic(t)x\psi_{x} -id(t)\psi + g(t)\psi  + h(t) (\left\vert \psi \right\vert ^{2} +  \left\vert \varphi \right\vert ^{2})\psi,
\end{eqnarray}
\begin{eqnarray}
\label{Eq1b}
i\varphi _{t} &=&a(t) \varphi_{xx} - b(t)x^2\varphi-ic(t)x\varphi_{x} -id(t)\varphi - g(t)\varphi  - h(t) (\left\vert \psi \right\vert ^{2} +  \left\vert \varphi \right\vert ^{2})\varphi.
\end{eqnarray}
Here, $\psi(x,t)$ and $\varphi(x,t)$ represent complex wavefunctions, and $a(t), b(t), c(t), d(t), g(t), h(t)$ are given time-dependent functions. The coefficient $a(t)$ represents the dispersion term and the opposite signs indicate that one component may experience normal dispersion, while the other experiences anomalous dispersion. The term $b(t)$ governs the time-dependent harmonic potential applied to the system. This configuration permits trapping of the $\psi$ component, whereas the opposite effect, an anti-trapping potential, in the $\varphi$ component promotes its spatial expansion. The third term, $c(t)$, governs the spatial dilation or compression of the solutions. The functions $d(t)$ and $g(t)$ represent, respectively, the gain or loss mechanism and an external linear potential. Ultimately, the function $h(t)$ controls the time-dependent nonlinear coupling between $\psi$ and $\varphi$, introducing an asymmetric interplay of self- and cross-phase modulation.

Assuming that $\psi$ and $\varphi$ solve the system (\ref{Eq1a})-(\ref{Eq1b}), we define the associated quantities below:
\begin{equation}
  N(t) =  \int_{-\infty}^{\infty}\left(|\psi(x,t)|^2 + |\varphi(x,t)|^2 \right) \ dx,  \quad \quad  (\mbox{mass or  $L^2$ norm}),
\end{equation}
and 
\begin{equation}
   L(t) =  \int_{-\infty}^{\infty} \operatorname{Im} \left( \psi^{*}(x,t) \, \psi_x(x,t) + \varphi^{*}(x,t) \, \varphi_x(x,t) \right) \, dx,  \quad \quad  (\mbox{linear momentum}), 
\end{equation}
 with the star denoting the complex conjugation. Then, formally, these solutions satisfy the following balance relations (assuming solutions that vanish as $|x| \rightarrow \infty$)
\begin{equation}
  N(t)  = \exp\left[\int_{0}^{t}\left(c(s)-2d(s)\right)\ ds\right] N(0), 
\end{equation}
and 
\begin{equation}
\frac{dL}{dt} = -2b(t) \int_{-\infty}^{\infty} x \left( |\psi|^2 - |\varphi|^2 \right) dx 
- 2d(t) L(t).
\end{equation}
It should be noted that the quantity $N(t)$ remains conserved when $c(t) = 2d(t)$, whereas the conservation of $L(t)$ requires that both $b(t)$ and $d(t)$ vanish. 

This work, grounded in the context mentioned above, is centered on analyzing the integrability and the intrinsic nonlinear behavior of system (\ref{Eq1a})-(\ref{Eq1b}). To begin our investigation into the system’s integrability, we employ bifurcation theory to identify a class of traveling-wave solutions. As will be shown later, these solutions naturally accommodate profiles expressed in terms of Jacobi elliptic functions, as well as bell-shaped and kink-type solitons. The presence of such solutions is conditioned by the behavior of the coefficients $g(t)$ and
$h(t)$. Secondly, the sensitivity analysis carried out in this study indicates that the obtained solutions could be relevant in both physical and numerical contexts. Nonetheless, the system may exhibit chaotic dynamics under the influence of external perturbations. Finally, modulational instability is explored in detail, initially considering the case $b(t)=0$, which permits a fully analytical treatment of the phenomenon. This preliminary analysis serves as a foundation for subsequent numerical simulations, which demonstrate that the system's dynamics under a harmonic potential still exhibit modulational instability, closely resembling the behavior observed in the non-harmonic scenario, at least during the early stages of the evolution. 

The presentation is structured as follows. In Section \ref{Sect2}, a comprehensive bifurcation analysis is performed to investigate the qualitative behavior of the solutions to system (\ref{Eq1a})-(\ref{Eq1b}). The explicit formulation of these solutions is described in Section \ref{Sect3}. Section \ref{Sect4} addresses the sensitivity analysis and examines the emergence of chaotic behavior when the coupled system is exposed to external perturbations. Section \ref{Sect5}  explores modulational instability, while the conclusions and final remarks of this investigation are presented in Section \ref{Sect6}. Section \ref{Sect7} serves as the appendix to this work.

\section{ Bifurcation Analysis of the CNLS Equations} \label{Sect2}

This section focuses on fully characterizing traveling wave solutions of (\ref{Eq1a})-(\ref{Eq1b}) by applying the bifurcation theory of an associated planar dynamical system. This bifurcation analysis reveals the existence of closed, homoclinic, and heteroclinic orbits, which are indicative of periodic solutions (Jacobi elliptic functions), bell-shaped, and kink-shaped solutions, respectively (the next section provides the explicit form of these solutions). 

\begin{figure}[h!]
\centering
\subfigure[Saddle point.]{\includegraphics[scale=0.65]{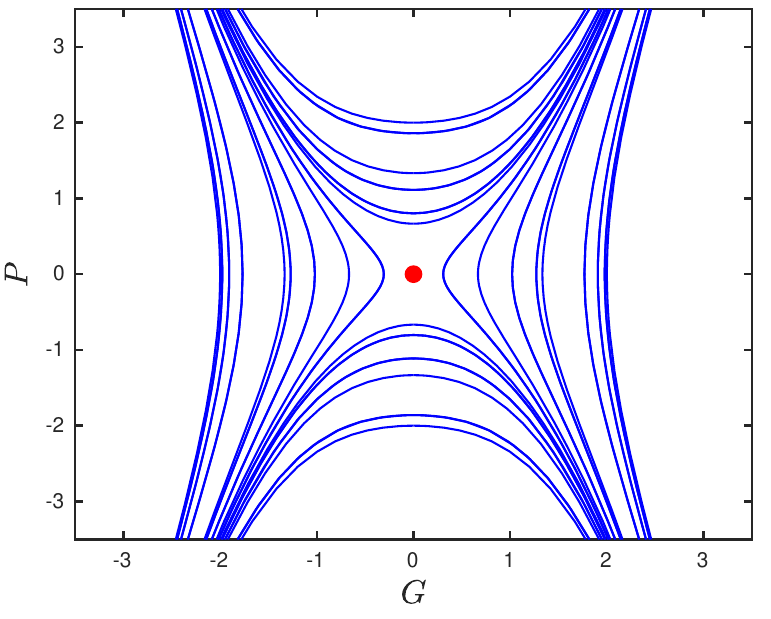}}
\subfigure[Center point.]{\includegraphics[scale=0.65]{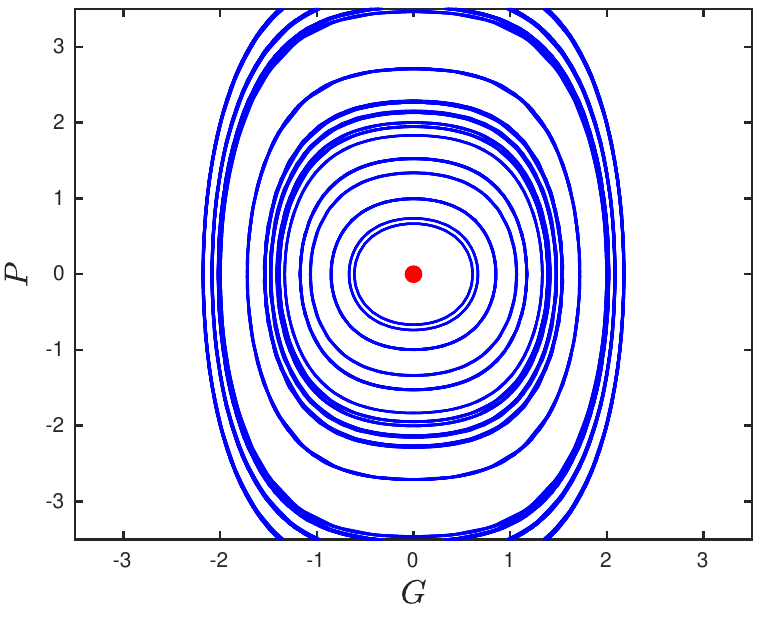}}
\caption{The phase portraits of system (\ref{PlanarSystem}) with (a) $g_0 = 1, h_0 = 0.5,$ and (b) $g_0 = -1, h_0 = -0.5.$ }\label{Fig1}
\end{figure}

\begin{figure}[h!]
\centering
\subfigure[Saddle/center points.]{\includegraphics[scale=0.65]{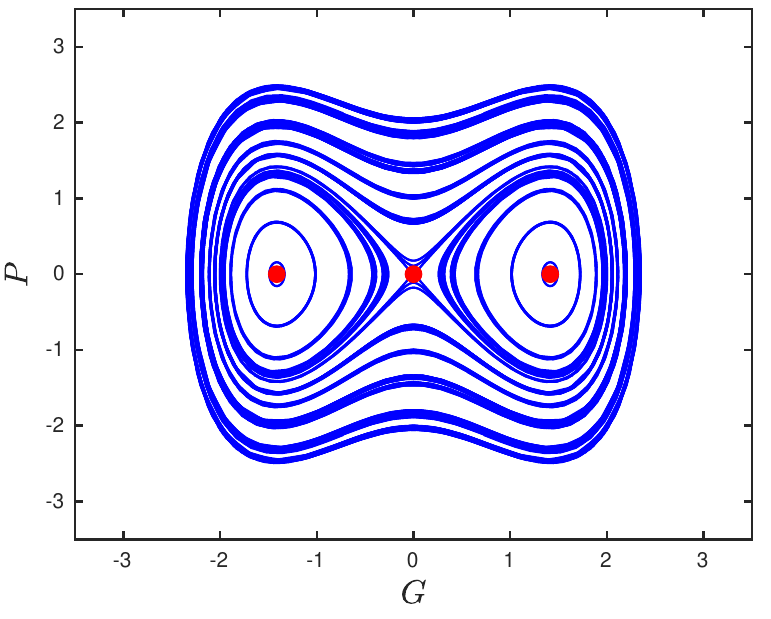}}
\subfigure[Center/saddle points.]{\includegraphics[scale=0.65]{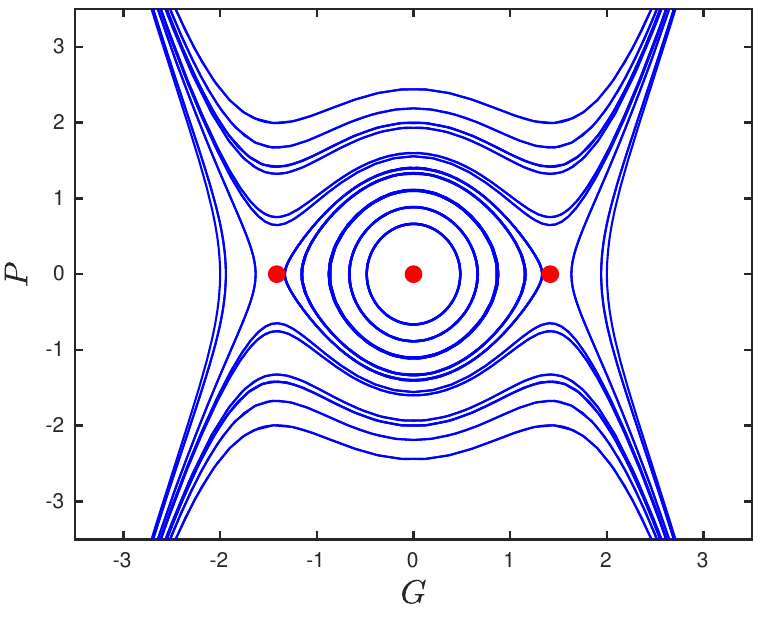}}
\caption{The phase portraits of system  (\ref{PlanarSystem}) with (a) $g_0 = 2, h_0 = -0.5,$ and (b) $g_0 = -2,  h_0 = 0.5.$ }\label{Fig2}
\end{figure}

Owing to the nature of the system (\ref{Eq1a})-(\ref{Eq1b}), let us consider the following transformations
\begin{equation}
  \psi(x,t) = e^{i\theta(x,t)}G\left(\xi\right),  \quad  \varphi(x,t) = e^{-i\theta(x,t)}G\left(\xi\right),    \label{Trans1}
\end{equation}
where $\theta(x,t) = \alpha(t)x^2 + \beta(t)x + \gamma(t) $ and $\xi = \beta(t)x + 2\gamma(t).$ Then by computing the derivatives of the functions $\psi(x,t)$ and $\varphi(x,t)$, we get
\begin{equation}
  \psi_x =  e^{i\theta(x,t)}\left( \beta G^{\prime} + 2i\alpha xG + i \beta G  \right),  \quad  \varphi_x =  e^{-i\theta(x,t)}\left( \beta G^{\prime} - 2i\alpha xG - i \beta G  \right), \label{Trans2}
\end{equation}
\begin{equation}
    \psi_{xx} = e^{i\theta(x,t)}\left( 4i \alpha \beta xG^{\prime} - 4\alpha^{2}x^2 G -4\alpha \beta xG + 2i\beta^{2}G^{\prime} - \beta^{2}G + \beta^2 G^{\prime \prime} + 2i \alpha G\right), \label{Trans3}
\end{equation}
\begin{equation}
    \varphi_{xx} = e^{-i\theta(x,t)}\left( -4i \alpha \beta xG^{\prime} - 4\alpha^{2}x^2 G -4\alpha \beta xG - 2i\beta^{2}G^{\prime} - \beta^{2}G + \beta^2 G^{\prime \prime} - 2i \alpha G\right), \label{Trans4}
\end{equation}
\begin{equation}
   \psi_t = e^{i\theta(x,t)}\left(i\dot{\alpha} x^2 G + i\dot{\beta}x G + i \dot{\gamma} G + \dot{\beta}x G^{\prime} + 2\dot{\gamma}G^{\prime} \right), \label{Trans5}
\end{equation}
\begin{equation}
    \varphi_t = e^{-i\theta(x,t)}\left(-i\dot{\alpha} x^2 G - i\dot{\beta}x G - i \dot{\gamma} G + \dot{\beta}x G^{\prime} + 2\dot{\gamma}G^{\prime} \right). \label{Trans6}
\end{equation}
Inserting the previous expressions (\ref{Trans1})-(\ref{Trans6}) into the equations (\ref{Eq1a})-(\ref{Eq1b}), we get the following result 
\begin{equation}
    G^{\prime \prime} = g_0 G + 2h_0 G^3, \label{EqG}
\end{equation}
provided $d(t) = -2a(t)\alpha(t)$, $g(t) = g_0 a(t)\beta^{2}(t)$, $h(t) = h_0 a(t)\beta^{2}(t)$ with $g_0, h_0 \not = 0$ are real constants and the following Riccati system is satisfied
\begin{equation}
   \dot{\alpha}(t) + b(t) + 2c(t)\alpha(t) +  4a(t)\alpha^{2}(t) = 0, \label{Ricati1}
\end{equation}
\begin{equation}
   \dot{\beta}(t) + (c(t) + 4a(t)\alpha(t))\beta(t)= 0, \label{Ricati2}
\end{equation}
\begin{equation}
    \dot{\gamma}(t) + a(t)\beta^2(t) = 0. \label{Ricati3}
\end{equation}
The assumption on the coefficient $d(t)$, implies that equations (\ref{Ricati1})-(\ref{Ricati2})  become
\begin{equation}
   \dot{\alpha}(t) + b(t) + 2(c(t) - d(t))\alpha(t) = 0, \label{Ricati1_1}
\end{equation}
\begin{equation}
   \dot{\beta}(t) + (c(t) -2d(t))\beta(t)= 0.
   \label{Ricati2_2}
\end{equation}
In these terms, the solution of the Riccati system is given by
\begin{equation}
   \alpha(t) = \alpha(0)e^{-E(t)}-e^{-E(t)}\int_{0}^{t}b(s) e^{E(s)} ds,  \quad E(t) = 2\int_{0}^{t}(c(s)-d(s))ds,
\end{equation}
\begin{equation}
   \beta(t) = \beta(0)e^{-\int_{0}^{t}(c(s)-2d(s))ds}, \quad \gamma(t) = \gamma(0) - \int_{0}^{t}a(s)\beta^{2}(s) ds.
\end{equation}
Now, by defining $P = G^{\prime}(\xi)$,  the equation (\ref{EqG}) can be rewritten as
  \begin{equation}
      \begin{cases} G^{\prime} = P, \\ P^{\prime} = g_0G + 2h_0 G^3. \end{cases}  \label{PlanarSystem}
  \end{equation}
  Therefore,
  \begin{itemize}
      \item[(i)]{If $g_0/h_0 > 0,$ the system (\ref{PlanarSystem}) admits the unique fixed point $(0,0).$}
      \item[(ii)]{If $g_0/h_0 < 0,$ the fixed points are $(0,0)$, $(\sqrt{-g_0/2h_0},0)$ and $(-\sqrt{-g_0/2h_0},0).$}
  \end{itemize}
The determinant of the Jacobian matrix is given by
\begin{equation}
    D(G,P) = \begin{vmatrix}
0 & 1 \\
g_0 + 6h_0G^2 & 0 
\end{vmatrix} = -g_0-6h_0G^2.
\end{equation}
In case $(i)$, we find that $(0,0)$ is a saddle point if $g_0,h_0>0,$ and a center point if $g_0,h_0<0.$ In case ($ii$), the fixed point $(0,0)$ is a saddle point if $g_0> 0, h_0<0$ and a center point if $g_0< 0, h_0>0.$ The last two fixed points $(\pm \sqrt{-g_0/2h_0},0)$  are center points if $g_0> 0, h_0<0$ and saddle points if $g_0<0, h_0>0.$  The phase portraits associated with cases $(i)$ and $(ii)$ of the dynamical system (\ref{PlanarSystem}) are illustrated in Figures  \ref{Fig1} and \ref{Fig2} (the fixed points of the system are shown as red dots). The bifurcation analysis reveals that the behavior of the dynamical system (\ref{PlanarSystem}), and consequently the nature of the solutions of (\ref{Eq1a})-(\ref{Eq1b}), is determined by the values of the external potential coefficient $g_0$ and the nonlinearity coefficient $h_0$. 

As a conclusion to this section, we find that the existence of closed orbits (Figure \ref{Fig1}(b) and Figure \ref{Fig2}), homoclinic orbits (Figure \ref{Fig2}(a)) and heteroclinic orbits (Figure \ref{Fig2}(b)) can be ensured by appropriate choices of the parameters $g_0$ and $h_0$. This will undoubtedly play a key role in the construction of explicit solutions in the next section.

\section{Explicit Solutions for the Coupled NLS System} \label{Sect3}
In the current section, the goal is to use the system's Hamiltonian structure (\ref{PlanarSystem}) to prove the existence (by the explicit construction) of Jacobi elliptic solutions. The discovery of such solutions will prove that the system has perfect intrinsic mechanisms for generating plane waves, modulated periodic waves, and solitons. 

To accomplish the aim, let us note that the system (\ref{PlanarSystem}) is equivalent to the following  Hamilton's equations 
\begin{equation}
   \frac{\partial H}{\partial P} = G^{\prime}, \quad \quad   \frac{\partial H}{\partial G} = -P^{\prime},
\end{equation}
where the Hamiltonian has the form
$$H(G,P) = \frac{1}{2}P^2 -\frac{g_0}{2}G^2 - \frac{h_0}{2}G^4.$$
Now, we consider the surface  
$$H(G,P) = \frac{1}{2}P^2 -\frac{g_0}{2}G^2 - \frac{h_0}{2}G^4 = H_0, \quad \quad (H_0 \ \mbox{is}  \ \mbox{a} \ \mbox{constant}),$$
and we will demonstrate how the various kinds of solutions rely on the values of $H_0.$ 

\begin{figure}[h!]
\centering
\subfigure[Profile of $|\psi|^2$.]{\includegraphics[scale=0.28]{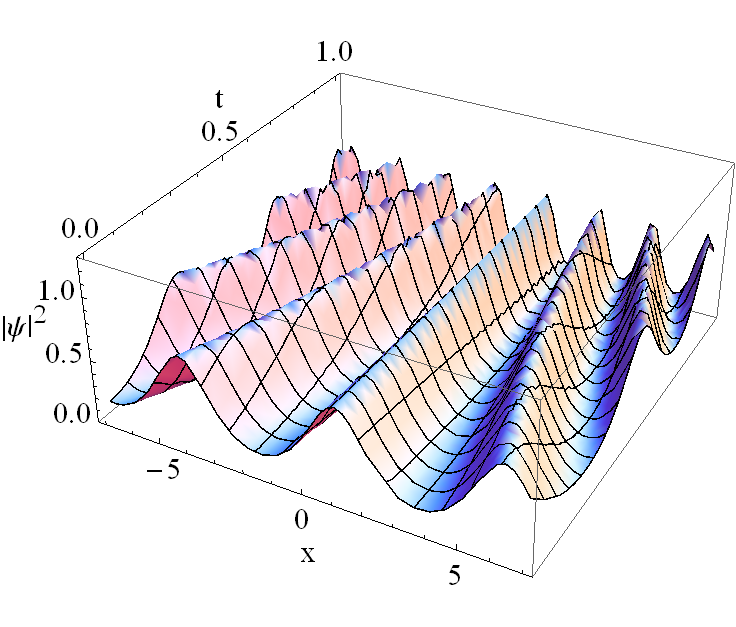}}
\subfigure[Profile of $|\psi|^2$.]{\includegraphics[scale=0.20]{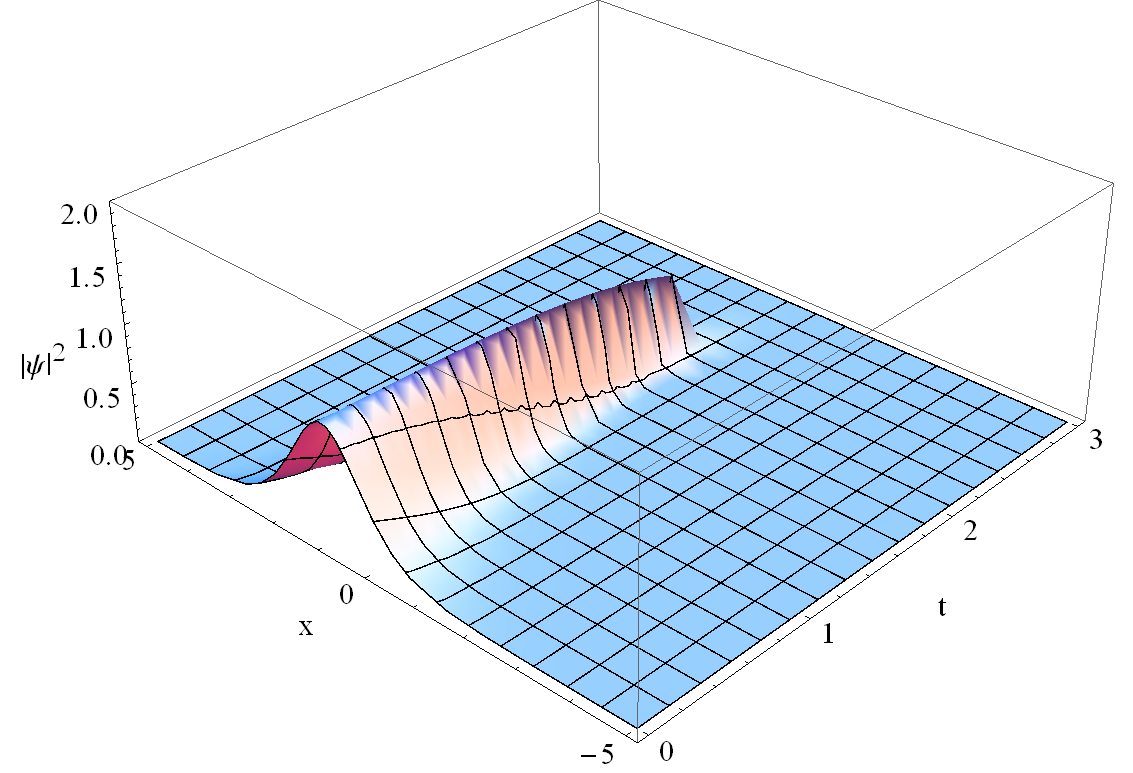}}
\subfigure[Profile of $|\psi|^2$.]{\includegraphics[scale=0.28]{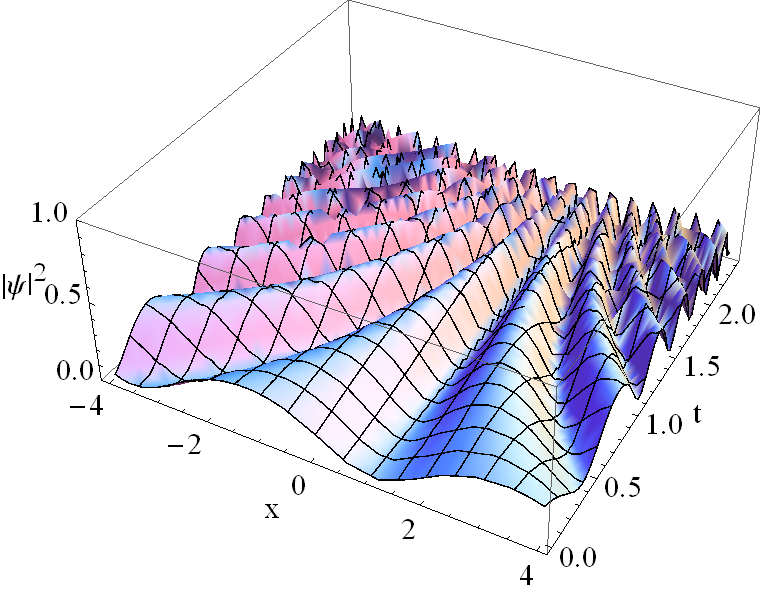}}
\subfigure[Profile of $|\psi|^2$.]{\includegraphics[scale=0.28]{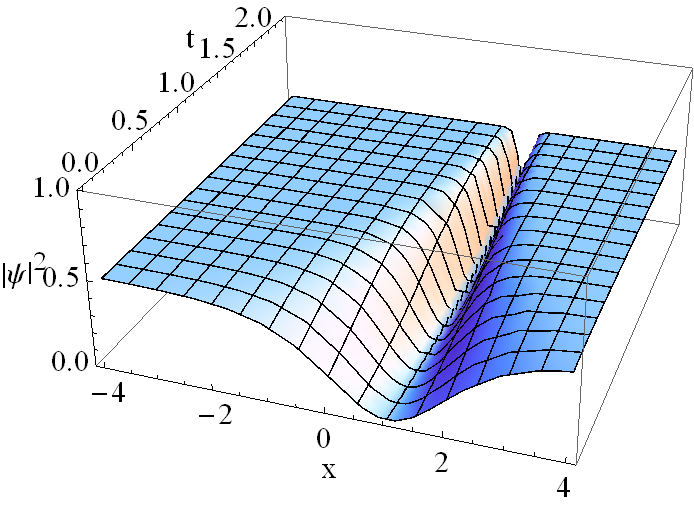}}
\caption{Solutions for the system (\ref{Eq1a})-(\ref{Eq1b}) with $a(t) =  0.5e^t,$   $b(t) = e^t,$ $c(t) = d(t) = 1,$ $g(t) = \frac{g_0}{2}e^t,$  $h(t) = \frac{h_0}{2}e^t,$ $\alpha(t) = -e^t,$ $\beta(t) = e^t,$ $\gamma(t) = -0.5e^t,$ and $\xi_0 = 0.$ (a) Delta amplitude solution (\ref{Sol1a}) with $h_0 = -1, g_0 = 1,$ $H_0 = -\frac{1}{10}.$ (b) Bell-shaped soliton solution (\ref{Sol2a}) with $h_0 = -1, g_0 = 1.$ (c) Sine amplitude solution (\ref{Sol5a}) with $h_0 = 1, g_0 = -1,$ $H_0 = \frac{1}{10}.$ (d) Dark soliton solution (\ref{Sol6a}) with $h_0 = 1, g_0 = -1,$ $H_0 = \frac{1}{8}.$}\label{Fig3}
\end{figure}

\begin{enumerate}
    \item[]{\textbf{Case 1: Consider $h_0 < 0$ and $g_0 > 0.$} \\
    $(a).$ Assume that $H_0 \in \left(\dfrac{g_0^2}{8h_0},0\right),$ then the Hamiltonian possesses four distinct real roots (corresponding to points on the closed orbits) and it can be expressed as: 
    \begin{equation}
        P^2 = g_0G^2 + h_0G^4 + 2H_0 =  -h_0\left(G^2 - G_{1H_0}^2\right)\left(G_{2H_0}^2 - G^2\right), \label{Hamilt2}
    \end{equation} 
    where
    \begin{equation}
        G_{1H_0}^2 = \frac{g_0 - \sqrt{g_0^2-8h_0H_0}}{-2h_0}, \quad \quad  G_{2H_0}^2 = \frac{g_0 + \sqrt{g_0^2-8h_0H_0}}{-2h_0}.
    \end{equation}
    Now, by remembering the definition of $P$, and integrating the equation (\ref{Hamilt2}),  we have
\begin{equation}
    \int_{-G_{2H_0}}^G \frac{ds}{\sqrt{\left(s^2 - G_{1H_0}^2\right)\left(G_{2H_0}^2 - s^2\right)}} = \pm \sqrt{-h_0}(\xi - \xi_0), \label{IntEllip1}
\end{equation}
    \begin{equation}
    \int_{G}^{G2H_0} \frac{ds}{\sqrt{\left(s^2 - G_{1H_0}^2\right)\left(G_{2H_0}^2 - s^2\right)}} = \mp \sqrt{-h_0}(\xi - \xi_0), 
\end{equation}
   where $\xi_0$ is a constant. The integral in (\ref{IntEllip1}) (and similarly the second integral) can be written as
   \begin{equation}
      \int_{-G_{2H_0}}^G \frac{ds}{\sqrt{\left(s^2 - G_{1H_0}^2\right)\left(G_{2H_0}^2 - s^2\right)}} = -\frac{1}{G_{2H_0}}\int_{1}^{-G/G_{2H_0}}\frac{ds}{\sqrt{(1-s^2)(s^2-\hat{l}^2)}},
   \end{equation}
   with $\hat{l}^2 = 1-l^2,$ and $l = \frac{\sqrt{G_{2H_0}^2 - G_{1H_0}^2}}{G_{2H_0}}.$ Therefore, by the definition of the \emph{delta amplitude function} $\dn(u;l)$, we obtain the explicit solutions for the coupled NLS equations
   \begin{equation}
    \psi(x,t) = \pm G_{2H_0}\dn\left( G_{2H_0}\sqrt{-h_0}(\beta(t)x + 2\gamma(t)-\xi_0); \frac{\sqrt{G_{2H_0}^2 - G_{1H_0}^2}}{G_{2H_0}}  \right)e^{i\theta(x,t)},  \label{Sol1a}
   \end{equation}
\begin{equation}
    \varphi(x,t) = \pm G_{2H_0}\dn\left( G_{2H_0}\sqrt{-h_0}(\beta(t)x + 2\gamma(t)-\xi_0); \frac{\sqrt{G_{2H_0}^2 - G_{1H_0}^2}}{G_{2H_0}}  \right)e^{-i\theta(x,t)}.  \label{Sol1b}
   \end{equation}
   $(b).$ In this case, we consider $H_0 = 0.$ Then, the roots of the Hamiltonian are $G_{1H_0}^2 = 0$ and $G_{2H_0}^2 = -g_0/h_0$ (corresponding to points on the homoclinic). Following the same reasoning of above, we find the bell-shaped solutions (using that $\dn(u,1) = \sech(u)$)
\begin{equation}
    \psi(x,t) = \pm \sqrt{-g_0/h_0}\sech\left( \sqrt{g_0}(\beta(t)x + 2\gamma(t)-\xi_0) \right)e^{i\theta(x,t)}, \label{Sol2a} 
   \end{equation}
\begin{equation}
    \varphi(x,t) = \pm \sqrt{-g_0/h_0}\sech\left( \sqrt{g_0}(\beta(t)x + 2\gamma(t)-\xi_0) \right)e^{-i\theta(x,t)}.  \label{Sol2b}
   \end{equation}
$(c).$ Consider that $H_0 = \dfrac{g_0^2}{8h_0}$. Under this condition $H(G,P)$ vanishes at the fixed points, i.e., $G_{1H_0}^2 = G_{2H_0}^2 = -g_0/2h_0.$ Therefore, the coupled NLS system (\ref{Eq1a})-(\ref{Eq1b}) admits the following plane wave solutions:
\begin{equation}
    \psi(x,t) = \pm \sqrt{-g_0/2h_0}e^{i\theta(x,t)} \label{Sol3a} 
   \end{equation}
\begin{equation}
    \varphi(x,t) = \pm \sqrt{-g_0/2h_0}e^{-i\theta(x,t)}.  \label{Sol3b}
   \end{equation}
$(d).$ Assume that $H_0>0$. In such case, we are considering the super periodic orbits of the system. Then, the Hamiltonian can be written as
\begin{equation}
        P^2 =  -h_0\left(G^2 + G_{3H_0}^2\right)\left(G_{2H_0}^2 - G^2\right), \label{Hamilt3}
    \end{equation} 
    where
    \begin{equation}
        G_{3H_0}^2 = \frac{-g_0 + \sqrt{g_0^2-8h_0H_0}}{-2h_0}.
    \end{equation}
Integrating the equation (\ref{Hamilt3}), we get
\begin{equation}
    \int_{G_{2H_0}}^G \frac{ds}{\sqrt{\left(s^2 + G_{3H_0}^2\right)\left(G_{2H_0}^2 - s^2\right)}} = \pm \sqrt{-h_0}(\xi - \xi_0). \label{IntEllip2}
\end{equation}
Note that the last integral is equivalent to
{\footnotesize
\begin{equation}
    \int_{G_{2H_0}}^G \frac{ds}{\sqrt{\left(s^2 + G_{3H_0}^2\right)\left(G_{2H_0}^2 - s^2\right)}} = \frac{1}{\sqrt{G_{2H_0}^2 + G_{3H_0}^2}} \int_{1}^{G/G_{2H_0}} \frac{ds}{\sqrt{\left(1-s^2\right)\left(\hat{l}^2 + l^2 s^2\right)}},
\end{equation}}
with $l^2 = \frac{G_{2H_0}^2}{G_{2H_0}^2 + G_{3H_0}^2}$ and $\hat{l}^2  = 1 -l^2.$ In these terms, we can construct the following Jacobian elliptic solutions
{\footnotesize
\begin{equation}
    \psi(x,t) = \pm G_{2H_0} \cn\left(\sqrt{-h_0}\sqrt{G_{2H_0}^2 + G_{3H_0}^2}(\beta(t)x + 2\gamma(t) - \xi_0); \frac{G_{2H_0}}{\sqrt{G_{2H_0}^2 + G_{3H_0}^2}}\right) e^{i\theta(x,t)}, \label{Sol4a} 
   \end{equation}
\begin{equation}
    \varphi(x,t) = \pm G_{2H_0} \cn\left(\sqrt{-h_0}\sqrt{G_{2H_0}^2 + G_{3H_0}^2}(\beta(t)x + 2\gamma(t) - \xi_0); \frac{G_{2H_0}}{\sqrt{G_{2H_0}^2 + G_{3H_0}^2}}\right) e^{-i\theta(x,t)}, \label{Sol4b}
   \end{equation}}
where $\cn(u;l)$ is the \emph{cosine amplitude function}.
}
 \item[]{\textbf{Case 2: Consider $h_0 > 0$ and $g_0 < 0.$} \\
 $(a).$ Firstly, assume that $H_0 \in \left(0,\dfrac{g_0^2}{8h_0}\right).$ Then, 
\begin{equation}
        P^2 =  h_0\left(G^2 - G_{1H_0}^2\right)\left( G^2-G_{2H_0}^2 \right), \label{Hamilt4}
    \end{equation} 
with 
\begin{equation}
        G_{1H_0}^2 = \frac{-g_0 - \sqrt{g_0^2-8h_0H_0}}{2h_0}, \quad \quad  G_{2H_0}^2 = \frac{-g_0 + \sqrt{g_0^2-8h_0H_0}}{2h_0}. 
    \end{equation}
Now, integrating the equation (\ref{Hamilt4}), we have
{\footnotesize
\begin{equation}
    \int_{0}^G \frac{ds}{\sqrt{\left(G_{1H_0}^2 -s^2\right)\left(G_{2H_0}^2 - s^2\right)}} = \frac{1}{G_{2H_0}}\int_{0}^{G/G_{1H_0}} \frac{ds}{\sqrt{\left(1-s^2\right)\left(1 - l^2 s^2\right)}} = \pm \sqrt{h_0}(\xi-\xi_0),
\end{equation}}
where $l = \frac{G_{1H_0}}{G_{2H_0}}.$ We  conclude that the system (\ref{Eq1a})-(\ref{Eq1b}) admits  solutions in terms of the \emph{sine amplitude function} $\sn(u;l)$, i.e.,
{\small
\begin{equation}
    \psi(x,t) = \pm G_{1H_0} \sn\left(\sqrt{h_0}G_{2H_0}(\beta(t)x + 2\gamma(t) - \xi_0); \frac{G_{1H_0}}{G_{2H_0}}\right) e^{i\theta(x,t)}, \label{Sol5a} 
   \end{equation}
\begin{equation}
    \varphi(x,t) = \pm G_{1H_0} \sn\left(\sqrt{h_0}G_{2H_0}(\beta(t)x + 2\gamma(t) - \xi_0); \frac{G_{1H_0}}{G_{2H_0}}\right) e^{-i\theta(x,t)}.\label{Sol5b}
   \end{equation}}
$(b).$ In this case, let us assume that $H_0 = \dfrac{g_0^2}{8h_0}.$ The roots of the Hamiltonian are $G_{1H_0}^2 = G_{2H_0}^2 = -\frac{g_0}{2h_0}$ (corresponding to points on the heteroclinic orbit). Therefore, we obtain the equality
\begin{equation}
    \int_{0}^G \frac{ds}{\sqrt{\left(\frac{g_0}{2h_0} +s^2\right)^2}} = \frac{1}{\sqrt{-g_0/2h_0}}\int_{0}^{\sqrt{-2h_0/g_0}G }\frac{ds}{\sqrt{\left(1-s^2\right)^2}} = \pm \sqrt{h_0}(\xi-\xi_0).
\end{equation}
Using the identity $\sn(u;1) = \tanh(u),$ we find the solutions
\begin{equation}
    \psi(x,t) = \pm \sqrt{-g_0/2h_0} \tanh\left(\sqrt{-g_0/2}(\beta(t)x + 2\gamma(t) - \xi_0)\right) e^{i\theta(x,t)}, \label{Sol6a} 
   \end{equation}
\begin{equation}
    \varphi(x,t) = \pm \sqrt{-g_0/2h_0} \tanh\left(\sqrt{-g_0/2}(\beta(t)x + 2\gamma(t) - \xi_0)\right) e^{-i\theta(x,t)}.\label{Sol6b}
   \end{equation}
}
\end{enumerate}
The profile of some solutions is shown in Figure \ref{Fig3}. Only $|\psi|^2$ is displayed, since the profile of the other component, $\varphi$, exhibits identical behavior. A noteworthy point, evident from the previously reported solutions, is that the amplitude and width of the solutions can be effectively controlled through the parameters $g_0$ and $h_0$. This characteristic could be particularly useful in certain physical applications, for example, where the external linear potential $g(t)$ could be used to manipulate the number of particles in a Bose-Einstein condensate.

\section{Sensitivity analysis and existence of chaotic behavior} \label{Sect4}
In this section, we explore the stability of the solutions of (\ref{Eq1a})-(\ref{Eq1b}) with respect to the initial conditions and external perturbations. This type of analysis is essential to ensure the physical and numerical application of the solutions reported in this work. 

\begin{figure}[h!]
\centering
\subfigure[$2D$ phase portrait without perturbation.]{\includegraphics[scale=0.7]{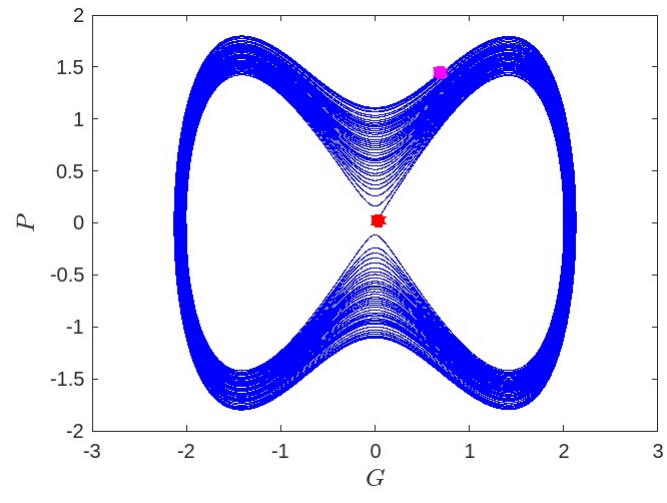}}
\subfigure[Dynamics of $G(\xi).$]{\includegraphics[scale=0.7]{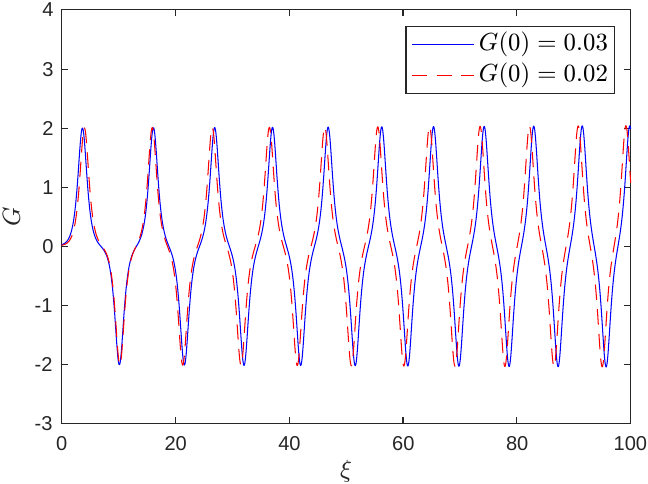}}
\subfigure[Dynamics of $P(\xi).$]{\includegraphics[scale=0.7]{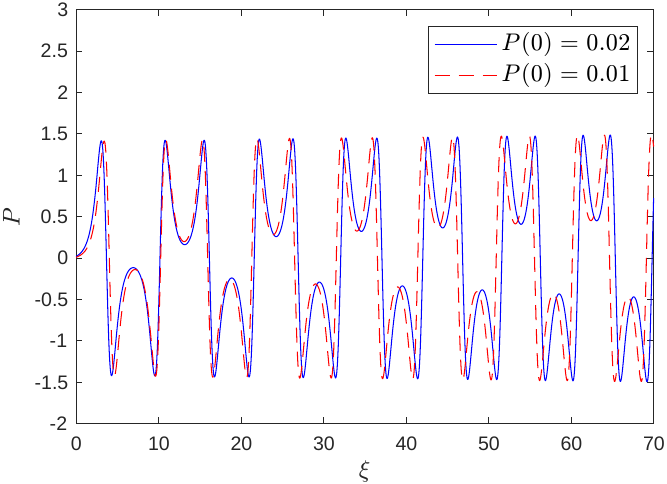}}
\subfigure[Spectrum of $G(\xi).$]{\includegraphics[scale=0.73]{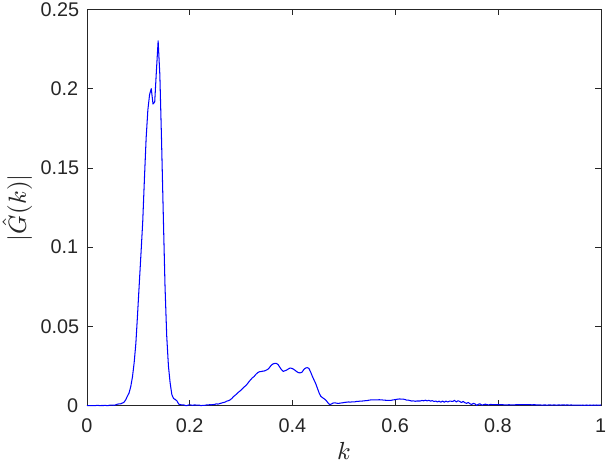}}
\caption{Dynamics of the unperturbed system (\ref{Perturbed_PlanarSystem}) ($\epsilon = 0$) with $h_0 = -0.5$ and $g_0 = 2$.  (a) $2D$ phase portrait of a system trajectory. The red and pink points correspond to the trajectory starting and ending positions, respectively. (b) Behavior of  $G(\xi)$  under two close initial conditions. (c) Same as (b) but for function $P(\xi)$. (d) Spectrum of $G(\xi)$ with conditions $G(0) = 0.03$ and $P(0) = 0.02$.  }\label{Fig4}
\end{figure}

In order to investigate the stability of the solutions against small fluctuations in the initial conditions, the system (\ref{PlanarSystem}) was numerically solved subject to the following values: $G(0) = 0.03$, $P(0) = 0.02$ and $G(0) = 0.02$, $P(0) = 0.01$. The results are shown in Figure \ref{Fig4}. The trajectory of the system with the values $G(0) = 0.03$ and $P(0) = 0.02$ is depicted in Figure \ref{Fig4}(a). Here, the red and pink dots indicate the starting and ending points of such a trajectory. The quasi-periodic behavior seen in the figure and, supported by its spectrum in (d), is in agreement with the previous bifurcation analysis. Figures \ref{Fig4}(b)-(c) exhibit the dynamics of $G$ and $P$, suggesting the robustness of the solutions when subjected to small changes in the initial conditions.

The second part of this section is devoted to examining the dynamics of the system under external perturbations. For that purpose, we explore the existence of chaos for the perturbed system 
\begin{equation}
      \begin{cases} G^{\prime} = P, \\ P^{\prime} = g_0G + 2h_0 G^3 + \epsilon \cos(K\xi), \end{cases}  \label{Perturbed_PlanarSystem}
  \end{equation}
 with $\epsilon$ representing the strength perturbation
 and $K$ the corresponding frequency (wavenumber). In this sense, system (\ref{Perturbed_PlanarSystem}) is equivalent to adding external perturbations $-\frac{\epsilon}{g_0}g(t) e^{i \theta}\cos(K \xi)$ and $\frac{\epsilon}{g_0}g(t) e^{-i \theta}\cos(K \xi)$ to equations (\ref{Eq1a}) and (\ref{Eq1b}), respectively. These types of perturbation correspond to a coherent space-time forcing in which certain frequency bands are excited over time \cite{Raza2023,Alfimov2002,Nail1997,Sulem1999,Tian2023,Wang2016}.
 
\begin{figure}[h!]
\centering
\subfigure[$2D$ phase portrait with perturbation.]{\includegraphics[scale=0.67]{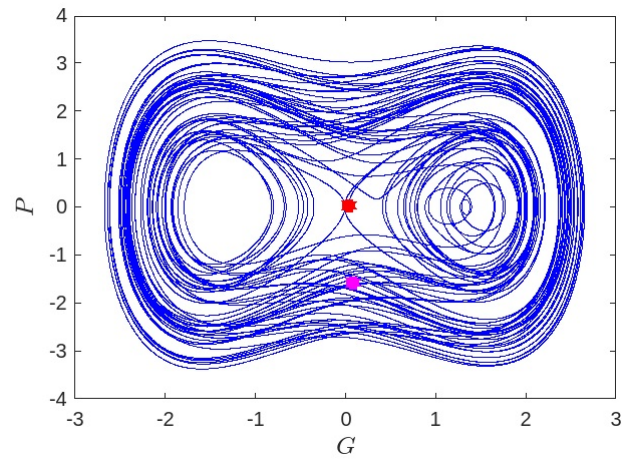}}
\subfigure[$3D$ phase portrait of the system.]{\includegraphics[scale=0.69]{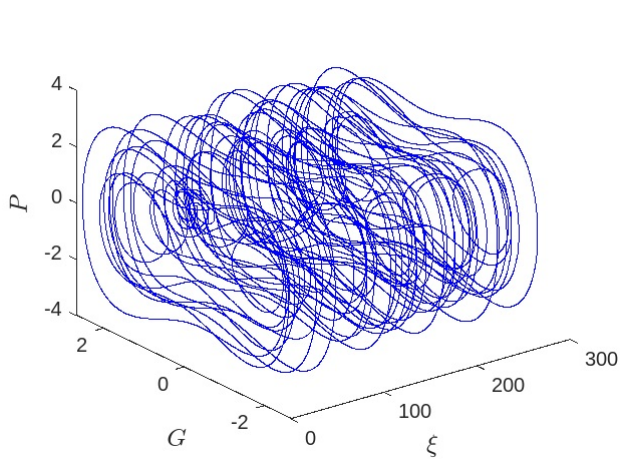}}
\subfigure[Dynamics of $G(\xi).$]{\includegraphics[scale=0.67]{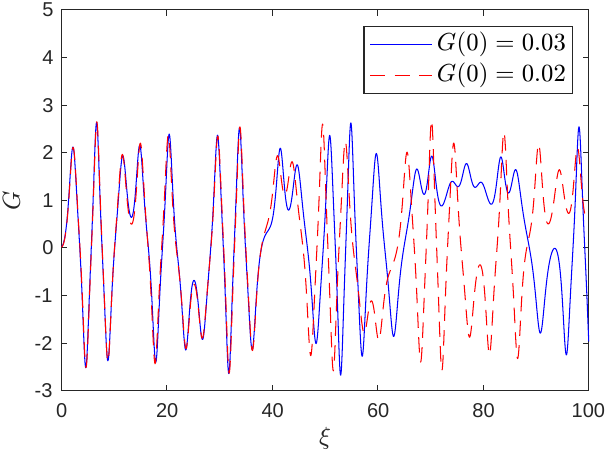}}
\subfigure[Dynamics of $P(\xi).$]{\includegraphics[scale=0.67]{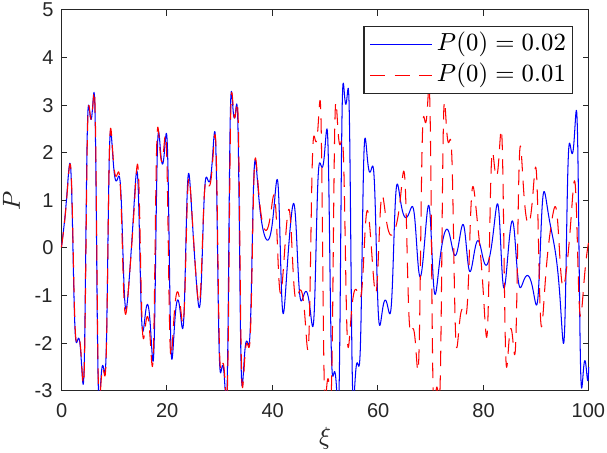}}
\subfigure[Spectrum of $G(\xi).$]{\includegraphics[scale=0.67]{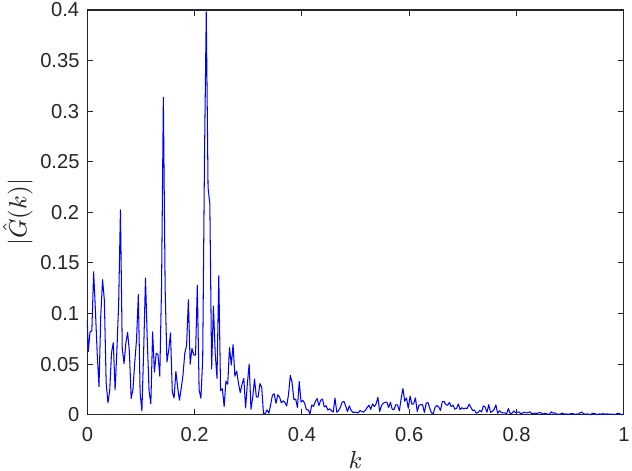}}
\caption{Dynamics of the system (\ref{Perturbed_PlanarSystem}) with $h_0 = -0.5$, $g_0 = 2$, and periodic perturbation with $\epsilon = 0.8$ and $K = 0.9$. (a) $2D$ phase portrait of the same trajectory shown in Figure 4(a). (b) $3D$ phase portrait of the  trajectory. (c) and (d) correspond to the dynamics of $G$ and $P$ for two close initial conditions. (e) Spectrum of $G(\xi)$ with conditions $G(0) = 0.03$ and $P(0) = 0.02$.    }\label{Fig5}
\end{figure}

 Figures \ref{Fig5}(a)-(b) exhibit a two-dimensional and three-dimensional visualization of the trajectory considered in Figure \ref{Fig4}(a)  when the system is exposed to a perturbation with $\epsilon = 0.8$ and $K = 0.9$. The erratic behavior of the trajectory demonstrates the presence of chaos in these specific settings. The numerical evidence presented in Figures \ref{Fig5}(c)-(d) indicates that the dynamics of the coupled NLS system vary significantly in response to small changes in the initial conditions a consequence of external perturbations. As is well known, this kind of instability is characteristic of chaotic systems (although the reciprocal is not true). Ultimately, the spectrum of the trajectory described in (a) is plotted. As shown in the figure, the spectrum of $G(\xi)$ is not concentrated around specific frequencies as in the quasi-periodic case (see Figure \ref{Fig4}(d)), and the graph is noisy but structured (presence of dominant frequencies and decays). These attributes are typical of chaotic systems.

As a consequence of the findings of this section, the constructed solutions have not only theoretical significance but also potential for physical and numerical applications. However, if the system is subject to external perturbations, chaotic behavior may occur, which destroys the dynamics of these solutions.  To mitigate these impacts, some type of control system must be implemented.

\section{Modulational instability}
\label{Sect5}
This part of the study aims to identify the conditions under which modulational instability arises in the system (\ref{Eq1a})-(\ref{Eq1b}). To this end, the initial analysis focuses on the case where $b(t)=0$ (no harmonic trap), providing insights that help elucidate the behavior in the more general case where 
$b(t) \not = 0$.

\subsection{Modulational instability for the CNLS without harmonic potential}
In this section, we investigate the modulation instability for the system (\ref{Eq1a})-(\ref{Eq1b}) without the quadratic potential, that is,
\begin{eqnarray}
\label{MODA}
i\psi _{t} &=&-a(t) \psi_{xx} -ic(t)x\psi_{x} -id(t)\psi + g(t)\psi  + h(t) (\left\vert \psi \right\vert ^{2} +  \left\vert \varphi \right\vert ^{2})\psi,
\end{eqnarray}
\begin{eqnarray}
\label{MODB}
i\varphi _{t} &=&a(t) \varphi_{xx} -ic(t)x\varphi_{x} -id(t)\varphi - g(t)\varphi  - h(t) (\left\vert \psi \right\vert ^{2} +  \left\vert \varphi \right\vert ^{2})\varphi.
\end{eqnarray}
For that purpose, we will consider the continuous wave (CW) solutions 
\begin{eqnarray}
  \psi_{0}(t) = A_0(t) e^{i\theta_1(t)}, \quad \quad  \varphi_{0}(t) = B_0(t) e^{i\theta_2(t)},   \label{CWSolutions}
\end{eqnarray}
with $\theta_1, \theta_2, A_0, B_0$ are real functions. By inserting the CW solutions into the system (\ref{MODA})-(\ref{MODB}), and equating the real and imaginary parts, we find
\begin{equation}
  A_0^{\prime}(t) = -d(t)A_0, \quad \quad   \theta_1^{\prime}(t) = -g(t)-h(t)\left(A_0^2 + B_0^2\right),
\end{equation}
\begin{equation}
  B_0^{\prime}(t) = -d(t)B_0, \quad \quad   \theta_2^{\prime}(t) = g(t)+h(t)\left(A_0^2 + B_0^2\right).
\end{equation}
Therefore, 
\begin{equation}
  \theta_1(t)  = \theta_{1}(0) - \int_{0}^{t} \left[g(s)+h(s)\left(A_{0}^{2}(s)  + B_{0}^{2}(s) \right)\right] ds, 
\end{equation}
\begin{equation}
  \theta_2(t)  = \theta_{2}(0) + \int_{0}^{t} \left[g(s)+h(s)\left(A_{0}^{2}(s)  + B_{0}^{2}(s) \right)\right] ds,
\end{equation}
\begin{equation}
 A_0(t) =  A_0(0)\exp\left(-\int_{0}^{t} d(s) \ ds\right), \quad    B_0(t) =  B_0(0)\exp\left(-\int_{0}^{t} d(s) \ ds\right).
\end{equation}
In these terms, the CW solutions (\ref{CWSolutions})  become
\begin{eqnarray}
  \psi_{0}(t) = A_0(0) \exp\left\{i\theta_{1}(0)-\int_{0}^{t} d(s) \ ds - i\int_{0}^{t} \left[g(s)+h(s)\left(A_{0}^{2}(s) + B_{0}^{2}(s)\right)\right] ds  \right\},  \label{CWa}
\end{eqnarray}
\begin{eqnarray}
  \varphi_{0}(t) = B_0(0) \exp\left\{i\theta_{2}(0)-\int_{0}^{t} d(s) \ ds + i\int_{0}^{t} \left[g(s)+h(s)\left(A_{0}^{2}(s) + B_{0}^{2}(s)\right)\right] ds  \right\}. \label{CWb}
\end{eqnarray}
As usual, we consider perturbations of the CW solutions as follows:
\begin{equation}
  \psi(x,t) = \psi_{0}(1 + \varepsilon(x, t)),  \label{Pert1} 
\end{equation}
\begin{equation}
  \varphi(x, t) = \varphi_{0}(1 + \delta(x, t)),  \label{Pert2}  
\end{equation}
with $\varepsilon(x,t)$, and $ \delta(x,t)$ are small perturbations. Inserting equations (\ref{Pert1})-(\ref{Pert2}) into the system (\ref{MODA})-(\ref{MODB}), and linearizing, we get 
\begin{equation}
i \varepsilon_{t} + a(t)\varepsilon_{xx} -h(t)\left(|\psi_0|^2 (\varepsilon + \varepsilon^{*}) + |\varphi_0|^2 (\delta + \delta^{*})\right)  = 0, \label{Li1} 
\end{equation}
\begin{equation}
i \delta_{t} - a(t) \delta_{xx} + h(t)\left(|\psi_0|^2 (\varepsilon + \varepsilon^{*}) + |\varphi_0|^2 (\delta + \delta^{*})\right)  = 0. \label{Li2} 
\end{equation}
Then, we look for solutions of (\ref{Li1})-(\ref{Li2}) in the form
\begin{equation}
  \varepsilon = r_1 e^{ikx-i\int_{0}^{t}w(s) ds} + r_2 e^{-ikx + i\int_{0}^{t}w(s) ds}, \quad \quad  \delta = r_3 e^{ikx-i\int_{0}^{t}w(s) ds} + r_4 e^{-ikx + i\int_{0}^{t}w(s) ds},
\end{equation}
where $k$ and $w$, are the wavenumber and the frequency, respectively. Here, $r_i$, are real constants that satisfy the homogeneous linear system $M \vec{r} = \vec{O}$ with
{\footnotesize
\begin{equation}
M = \left(
\begin{matrix}
w-a(t)k^2 -h(t)|\psi_0|^2 & -h(t)|\psi_0|^2 & -h(t)|\varphi_0|^2 & -h(t)|\varphi_0|^2\\
-h(t)|\psi_0|^2 & -w-a(t)k^2 -h(t)|\psi_0|^2 & -h(t)|\varphi_0|^2 & -h(t)|\varphi_0|^2\\
h(t)|\psi_0|^2 & h(t)|\psi_0|^2 & w+a(t)k^2 +h(t)|\varphi_0|^2 & h(t)|\varphi_0|^2 \\
h(t)|\psi_0|^2 & h(t)|\psi_0|^2 & h(t)|\varphi_0|^2 & -w+a(t)k^2 +h(t)|\varphi_0|^2 
\end{matrix}
\right),
\end{equation}
}
and $\vec{r} = (r_1,r_2,r_3,r_4)$. In order to admit no trivial solution, the determinant of $M$ should be zero, obtaining the following dispersion relation 
\begin{equation}
     \omega^2 =  a(t)^2k^4 +  a(t)h(t) k^2 \left(|\psi_0|^2 + |\varphi_0|^2\right) \pm |a(t)h(t)| k^2 \left(|\psi_0|^2 + |\varphi_0|^2\right). \label{DR}
\end{equation}
According to (\ref{DR}), some $\omega$ will be complex, and then the CW solutions are unstable, if the wave numbers $k$ satisfy the inequality
\begin{equation}
k^2 < -\frac{2h(t)}{a(t)}\left(|\psi_0|^2 + |\varphi_0|^2\right) = -\frac{2h(t)}{a(t)}\left(A_{0}^{2}(0) + B_{0}^{2}(0)\right)e^{-2\int_{0}^{t} d(s) \ ds}, \label{UnstRegion}
\end{equation}
provided $a(t)h(t)<0.$ Equation (\ref{UnstRegion}) reveals that modulation stability results from the interaction of the dispersion term, nonlinearity, dissipation and the initial amplitudes of the CW solutions. It should be noted that if the quotient between dispersion and nonlinearity terms is constant, then the region of unstable wavenumbers expands (contracts) as long as $d(t) <0$ ($d(t) > 0$). On the other hand, if the quotient of these terms is proportional to $e^{2\int_{0}^{t} d(s) \ ds}$, the region (\ref{UnstRegion}) is fixed. The gain parameter, $\Lambda$, is defined as
\begin{equation}
  \Lambda(k) = 2\mbox{Im}(w) =  2\sqrt{-a(t)^2k^4 -2  a(t)h(t) k^2 \left(A_{0}^{2}(0) + B_{0}^{2}(0)\right)e^{-2\int_{0}^{t} d(s) \ ds}}.  \label{Gain}
\end{equation}
For the previous relation, $\Lambda(k)$ is zero for $k = 0,\pm \sqrt{-\frac{2h(t)}{a(t)}\left(A_{0}^{2}(0) + B_{0}^{2}(0)\right)}e^{-\int_{0}^{t} d(s) \ ds}$, and the maximum value $\Lambda_{\max} = 2|h(t)|\left(A_{0}^{2}(0) + B_{0}^{2}(0)\right)e^{-2\int_{0}^{t} d(s) \ ds}$ is obtained at the wavenumbers $k_{\max} = \pm \sqrt{-\frac{h(t)}{a(t)}\left(A_{0}^{2}(0) + B_{0}^{2}(0)\right)}e^{-\int_{0}^{t} d(s) \ ds}$. In the subsections that follow, we will look at the dynamics of the gain spectrum, $\Lambda(k)$, when the dispersion and nonlinearity coefficients are periodic and exponential. In what follows $\theta_1(0) = \theta_2(0) = 0,$ and $A_0(0) = B_0(0) = 1.$

\subsubsection{Case $a(t) = 1 + \cos t$, $c(t) = 0,$ $d(t) = \pm t,$  $g(t) = \cos t$, and $h(t) = -2 - 2\cos t$} \label{Sect511} 
\begin{figure}[h!]
\centering
\subfigure[Dynamics for the gain $\Lambda^{+}$ given in (\ref{GainPer}). ]{\includegraphics[scale=0.52]{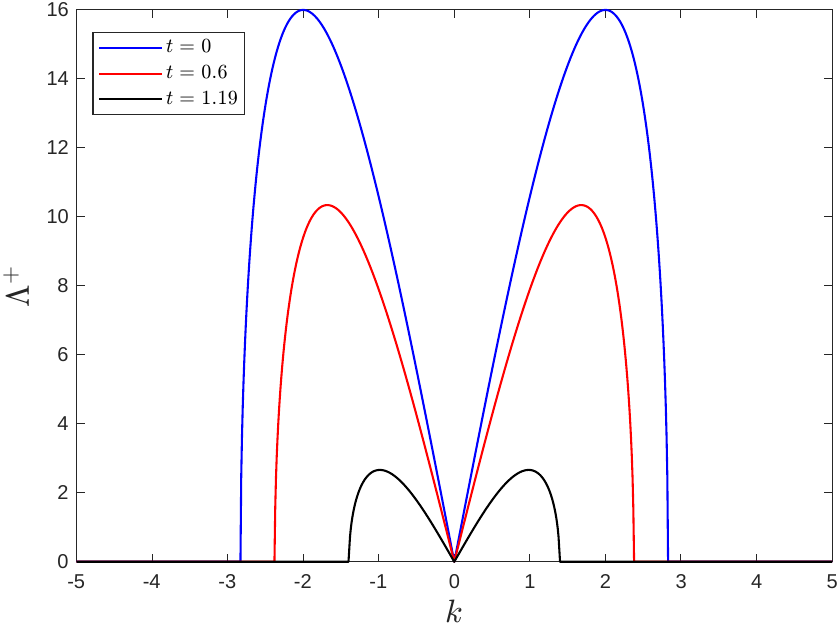}}
\subfigure[Dynamics for the gain $\Lambda^{-}$ given in (\ref{GainPer}).]{\includegraphics[scale=0.52]{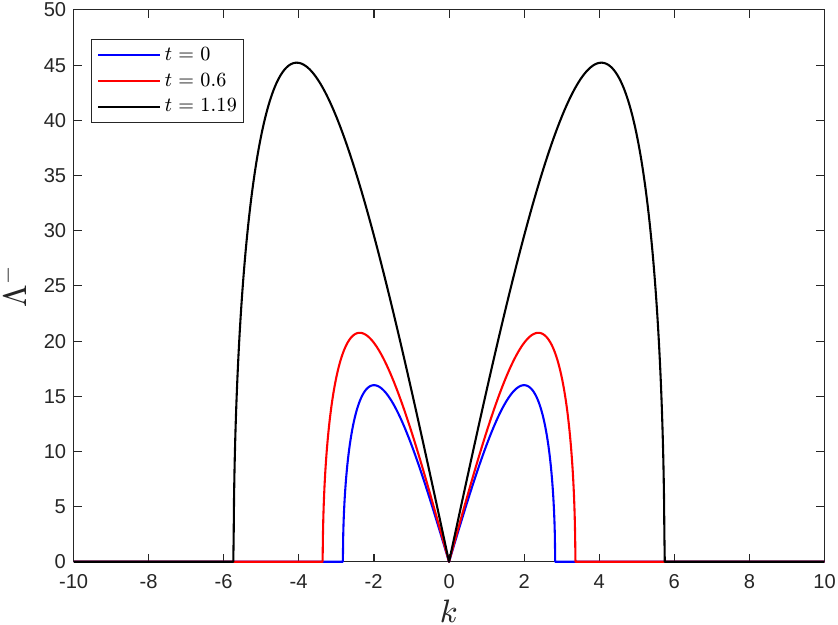}}
\subfigure[Dynamics for the gain $\Lambda$ given in (\ref{Gainexp}).]{\includegraphics[scale=0.54]{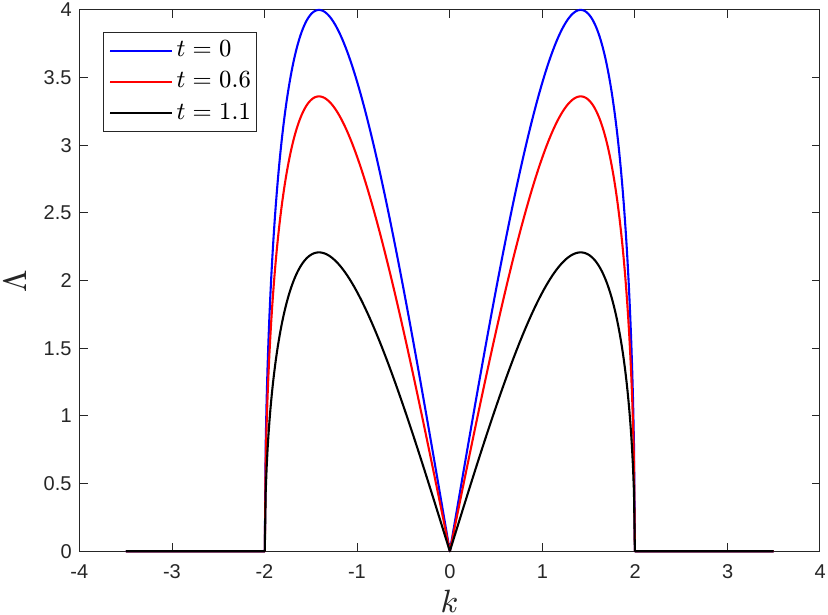}}
\subfigure[Dynamics for the gain $\Lambda$ given in (\ref{Gainexp2}). ]{\includegraphics[scale=0.57]{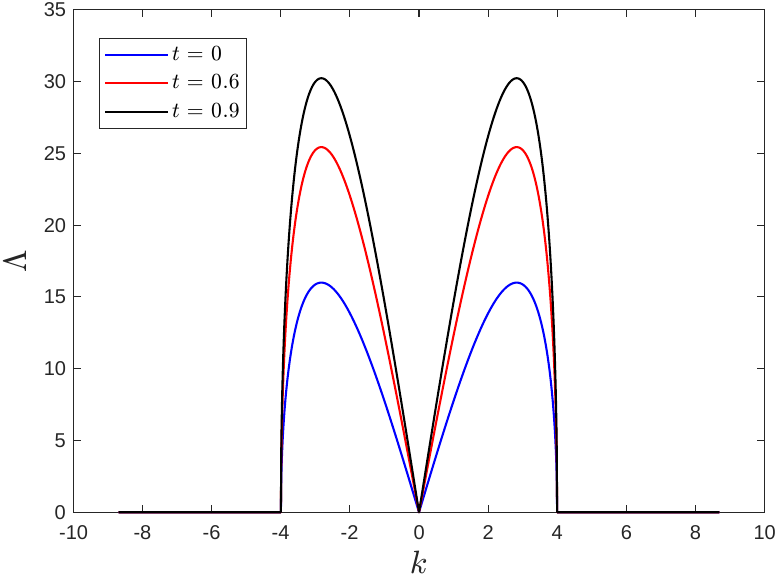}}
\caption{(a) In the long term, the positivity of $d(t) = t$ reduces both the amplitude of $\Lambda$ and the size of the instability region. (b) Contrary to the previous case, the amplitude of the gain and the the size of the instability region increase due to the influence of $ d(t) = -t$. (c) The interval of  unstable wavenumbers are fixed, but the amplitude of the gain decreases with time. (d) The interval of  unstable wavenumbers are fixed, but the amplitude of the gain increases with time.}
\label{Fig6}
\end{figure}

The selection of these coefficients ensures the existence of unstable wavenumbers. Here, we examine the dynamics of $\Lambda$ in relation to the gain/loss term $d(t)$. Based on the selection of $d(t)$, we obtain the two gain parameters $\Lambda^+$ and $\Lambda^-$:   
\begin{equation}
  \Lambda^{+} = 2|k|(1 + \cos t)\sqrt{8e^{-t^2}-k^2}, \quad \quad  \Lambda^{-} = 2|k|(1 + \cos t)\sqrt{8e^{t^2}-k^2} . \label{GainPer}
\end{equation}
The profiles of the gain parameters $\Lambda^+$ and $\Lambda^-$ are shown in Figures \ref{Fig6}(a)-(b). As can be seen, the function $d(t)$ can be used to control the instability of the system. A positive coefficient $d(t)$ prevents the instability (by mitigating unstable wavenumbers), whereas a negative coefficient $d(t)$ quite accelerates it(by unstably exciting new Fourier modes).     

\subsubsection{Case $a(t) = e^{-\frac{t^2}{2}}$, $c(t) = 0,$ $d(t) = t,$  $g(t) = \cos t$, and $h(t) = -e^{\frac{t^2}{2}}$} 
In this part, we demonstrate that, in certain scenarios, the coefficient $d(t)$ simply affects the amplitude of the gain spectrum and not the instability interval. The gain spectrum associated with these coefficients has the following form
\begin{equation}
  \Lambda = 2|k|e^{-\frac{t^2}{2}}\sqrt{4-k^2} .\label{Gainexp}
\end{equation}
The instability interval corresponds to wavenumbers with $|k| < 2$, while the amplitude decays exponentially over time, see Figure \ref{Fig6}(c).
\subsubsection{Case $a(t) = t + 1$, $c(t) = 0,$ $d(t) = -t,$  $g(t) = \cos t$, and $h(t) = -4(t + 1)e^{-t^2}$} \label{Sect513}
The gain spectrum for these variable coefficients is given by the expression
\begin{equation}
  \Lambda = 2|k|(t+1)\sqrt{16-k^2} .\label{Gainexp2}
\end{equation}
The form of $\Lambda$ reveals a linear growth of the amplitude in the variable $t$, while, as in the previous case, the region of unstable wavenumbers remains fixed ($|k|<4$). In the panel (d) of Figure \ref{Fig6}, we have plotted the profile of $\Lambda$. 

In order to verify the modulation instability characteristics pointed out beforehand, we have performed numerical simulations of the system (\ref{MODA})-(\ref{MODB}) with initial conditions 
\begin{equation}
  \psi(x,0) = 1 + \epsilon_1 \cos(kx), \quad \varphi(x,0) = 1 + \epsilon_2 \cos(kx),
\end{equation}

\begin{figure}[h!]
\centering
\subfigure[Profile of the function $|\psi|^2$.]{\includegraphics[scale=0.52]{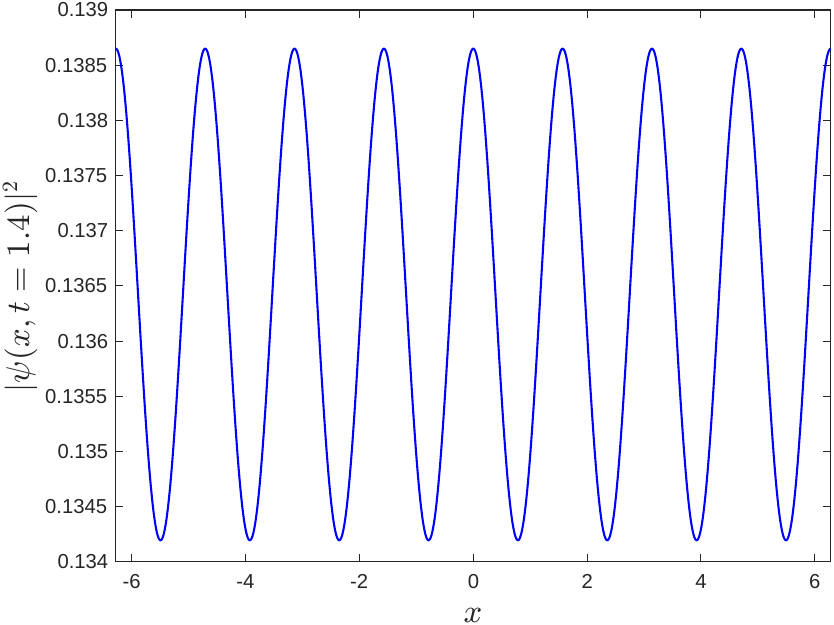}}
\subfigure[Profile of the function $|\psi|^2$ .]{\includegraphics[scale=0.52]{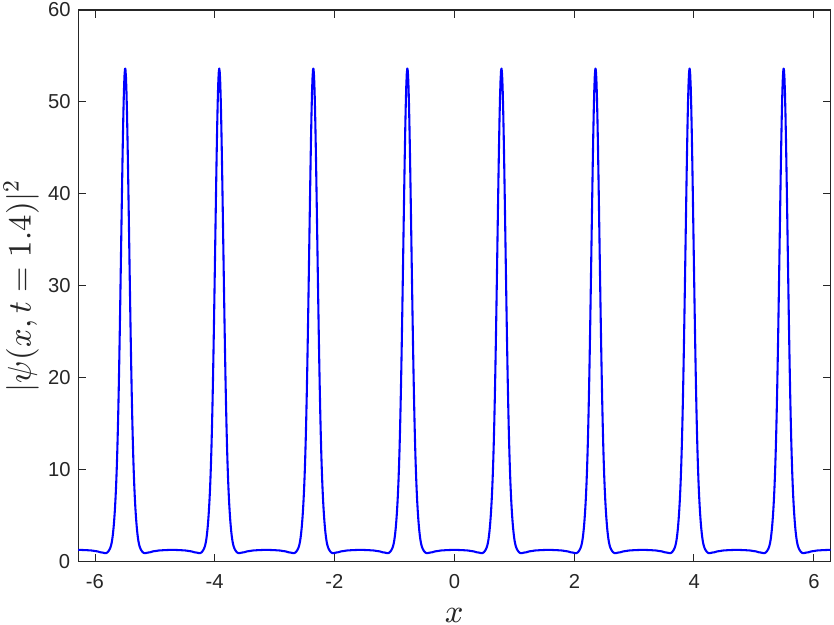}}
\caption{Numerical simulations for the system (\ref{MODA})-(\ref{MODB}) with $a(t) = 1 + \cos t$, $c(t) = 0,$  $g(t) = \cos t$, $h(t) = -2 - 2\cos t$  and initial conditions $\psi(x,0) = 1 + 10^{-2} \cos(4x) $ and $ \varphi(x,0) = 1 +  10^{-4}\cos(4x)$ (stable wavenumber). The profiles correspond to the time $t = 1.4$. (a) Influence of the function $d(t) = t$. (b) Influence of the function $d(t) = -t$.}\label{Fig7}
\end{figure}
\noindent  where $k$ and $\epsilon_i$ ($i = 1,2$) denote the wavenumber and amplitude of the perturbations, respectively. The numerical results for such a system with the same coefficients given in Section \ref{Sect511} are shown in Figure \ref{Fig7} (similar behavior is found for the function $\varphi$). The system undergoes an instability for wavenumbers that satisfy $|k| < \sqrt{8}$ (at time $t = 0$). For a wavenumber $ k = 4$ (stable), the coupled NLS responds to the dissipation term $d(t) = t$ with a small amplitude modulation, as shown in panel \ref{Fig7}(a). This is in agreement with the fact that a positive coefficient $d(t)$ tends to contract the instability region, so the initial wavenumber is never excited. In the other scenario, with $d(t) = -t$, the system (\ref{MODA})-(\ref{MODB}) admits the formation of a soliton wave train as a consequence of the expansion of the instability region; see Figure \ref{Fig7}(b).

\begin{figure}[h!]
\centering
\subfigure[Profiles of the functions $|\psi|^2$ and $|\varphi|^2$ .]{\includegraphics[scale=0.73]{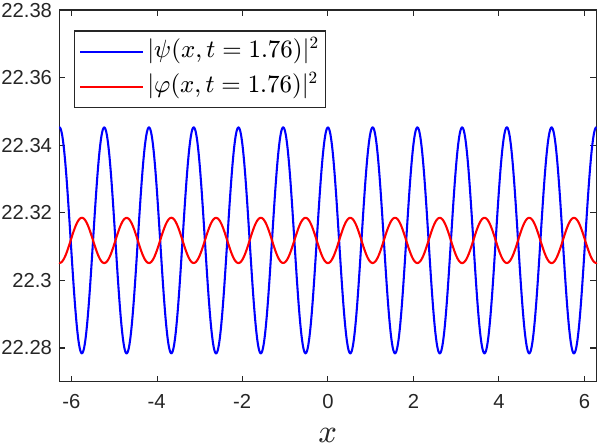}}
\subfigure[Profile of the function $|\psi|^2$.]{\includegraphics[scale=0.73]{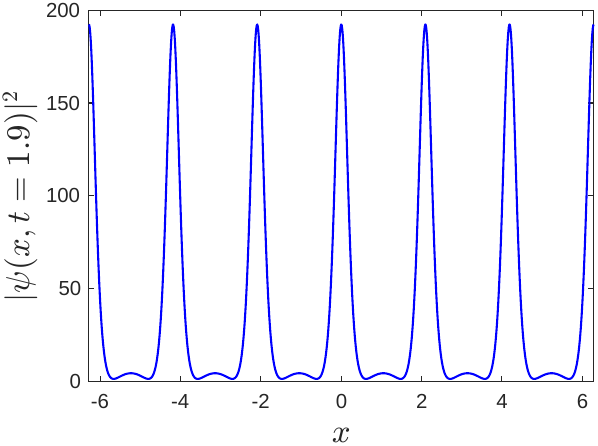}}
\caption{Numerical simulations for the system (\ref{MODA})-(\ref{MODB}) with $a(t) = 1 + t$, $c(t) = 0,$ $d(t) = -t,$  $g(t) = \cos t$, $h(t) = -4(1 + t)e^{-t^2}$  and initial conditions $\psi(x,0) = 1 + 10^{-3} \cos(kx) $ and $ \varphi(x,0) = 1 +  10^{-4}\cos(kx)$ . (a) Dynamics for the stable wavenumber $k = 6$ at $t = 1.76$. (b) Dynamics for the unstable wavenumber $k = 3$ at $t = 1.9$.}
\label{Fig8}
\end{figure}

The numerical simulations corresponding to the coefficients used in Section \ref{Sect513} can be found in Figure \ref{Fig8}. As we mentioned earlier, the system undergoes MI in a fixed interval $|k|< 4$. Figure \ref{Fig8} depicts the dynamics of systems perturbed by stable ($k = 6$) and unstable ($k = 3$) wavenumbers, as shown in panels (a) and (b), respectively. In the first part, the wavenumber remains stable, but the system is supplied with energy through the gain term $d(t) = -t$, and therefore a large amplitude modulation is observed. Finally, in the second part,  a soliton wave train is obtained as a result of the energy supply and the instability.

\subsection{Modulational instability for the   CNLS with harmonic potential}

This section of the article focuses on the MI for a coupled NLS system with a non-zero coefficient $b(t)$. Since this general model does not admit CW (time-uniform) solutions, we cannot perform the classical MI analysis. As a result, the plan is to connect the harmonic system ($b(t) \not = 0$) to the ones established in (5.1)-(5.2), in the hope that the analysis conducted with the former will reveal the presence of modulation instability in the latter.  The following theorem reveals a connection between both systems.

\begin{theorem} \label{Theorem1}
Suppose that $a(t), b(t), c(t), d(t), g(t)$ and $h(t)$ are time-dependent functions. Then, the general coupled NLS system
\begin{eqnarray}
\label{EqA}
i\psi _{t} &=&-a(t) \psi_{xx} + b(t)x^2\psi-ic(t)x\psi_{x} -id(t)\psi + g(t)\psi  + h(t) (\left\vert \psi \right\vert ^{2} +  \left\vert \varphi \right\vert ^{2})\psi,
\end{eqnarray}
\begin{eqnarray}
\label{EqB}
i\varphi _{t} &=&a(t) \varphi_{xx} - b(t)x^2\varphi-ic(t)x\varphi_{x} -id(t)\varphi - g(t)\varphi  - h(t) (\left\vert \psi \right\vert ^{2} +  \left\vert \varphi \right\vert ^{2})\varphi,
\end{eqnarray}
can be transformed into the constant coefficient system 
\begin{eqnarray}
\label{Moda}
iu _{\tau} + u_{\xi \xi} + id_0 u  - g_0 u  - h_0 (\left\vert u \right\vert ^{2} +  \left\vert v \right\vert ^{2})u = 0,
\end{eqnarray}
\begin{eqnarray}
\label{Modb}
iv _{\tau} - v_{\xi \xi} + id_0 v  + g_0 v  + h_0 (\left\vert u \right\vert ^{2} +  \left\vert v \right\vert ^{2})v = 0,
\end{eqnarray}
where $d_0, g_0, h_0$ are real constants and $\xi = \beta(t)x$ and $\tau = \gamma(t), $ are rescaled space and time variables.
\end{theorem}
\begin{proof}
Consider the substitutions 
\begin{equation}
   \psi(x,t) = \frac{1}{\sqrt{\mu(t)}}e^{i \alpha(t)x^2}u(\xi, \tau), \label{Transf1}
\end{equation}
\begin{equation}
   \varphi(x,t) = \frac{1}{\sqrt{\mu(t)}}e^{-i \alpha(t)x^2}v(\xi, \tau). \label{Transf2}
\end{equation}
Then, the first derivatives of the function $\psi$ are given by
\begin{equation}
  \psi_x = \mu^{-1/2}(t)e^{i\alpha(t)x^2} \left[ \beta u_\xi + 2\alpha i x u  \right], 
\end{equation}
\begin{equation}
  \psi_{xx} = \mu^{-1/2}(t)e^{i\alpha(t)x^2} \left[ 4\alpha \beta i x u_\xi -4\alpha^2  x^2 u + \beta^2 u_{\xi \xi} + 2 \alpha iu \right],
\end{equation}
\begin{equation}
\psi_t = \mu^{-1/2}(t)e^{i\alpha(t)x^2} \left[-\frac{1}{2}\frac{\dot{\mu}}{\mu}u + i \dot{\alpha}x^2 u + \dot{\beta} x u_\xi + \dot{\gamma}u_\tau \right].
\end{equation}
Inserting these expressions in the equation (\ref{EqA}), we have
\begin{align}
    x^2 u \left(-\dot{\alpha}  - 4a \alpha^2  -b -2c\alpha \right) &+ ix u_{\xi}\left( \dot{\beta} + 4a \alpha \beta  + c \beta \right) +iu \left(2a \alpha - \frac{1}{2}\frac{\dot{\mu}}{\mu}\right) \\
    &+  \left( i \dot{\gamma} u_\tau + a\beta^2 u_{\xi \xi} + idu -gu - h\mu^{-1}(|u|^2 + |v|^2)u  \right) = 0.
\end{align}
Therefore, the equation (\ref{Moda}) is obtained provided the following Riccati system holds 
\begin{equation}
   \dot{\alpha}(t) + b(t) + 2c(t)\alpha(t) +  4a(t)\alpha^{2}(t) = 0, \label{Riccati1}
\end{equation}
\begin{equation}
   \dot{\beta}(t) + (c(t) + 4a(t)\alpha(t))\beta(t)= 0, \label{Riccati2}
\end{equation}
\begin{equation}
    \dot{\gamma}(t) - a(t)\beta^2(t) = 0, \label{Riccati3}
\end{equation}
\begin{equation}
  \dot{\mu}(t) -4a(t)\alpha(t) \mu(t)= 0,  \label{Riccati4} 
\end{equation}
and $d(t) = d_0 a(t) \beta^2(t)$, $g(t) = g_0 a(t) \beta^2(t) $ and $h(t) = h_0 a(t) \beta^2(t) \mu(t)$. Subject to the same system of equations and similar reasoning, equation (\ref{Modb})  can be deduced. We can combine the equations (\ref{Riccati1}) and (\ref{Riccati4}) to produce the second order ODE
\begin{equation}
   \ddot{\mu} + \left(2c-\frac{\dot{a}}{a}\right)\dot{\mu} + 4ab\mu = 0.
\end{equation}
This last equation is named in the literature as the \emph{characteristic equation}. The characteristic equation can be used to solve the Riccati system, as shown in the appendix of this paper.  
\end{proof}

\begin{figure}[h!]
\centering
\subfigure[Profile of the function  $|\psi|^2$.]{\includegraphics[scale=0.8]{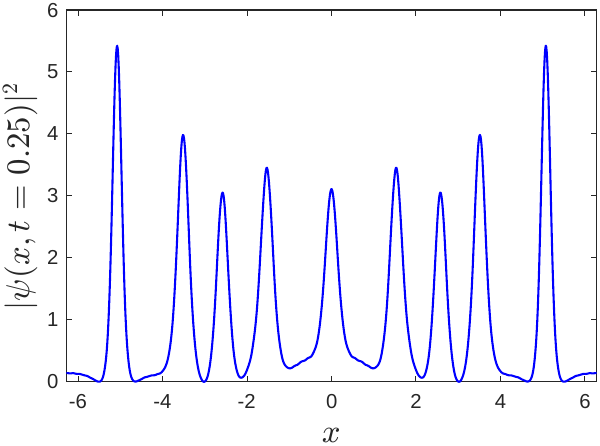}}
\subfigure[Profile of the function  $|\psi|^2$.]{\includegraphics[scale=0.8]{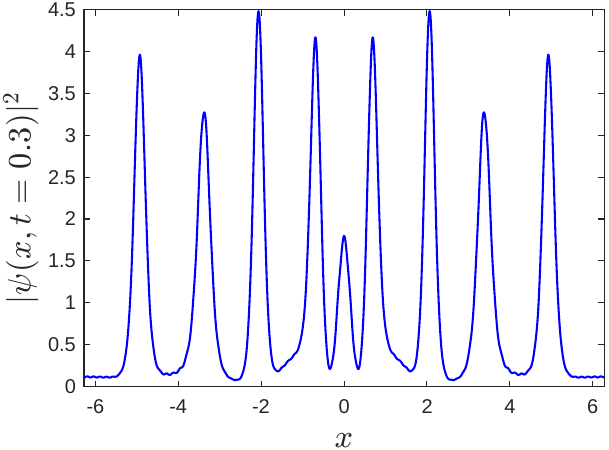}}
\caption{The snapshots of squared modulus of the solution for $\epsilon = 10^{-4}$ and $d_0 = -10^{-2}$. (a) Dynamics for the  stable wavenumber $k = 4$ at the instant $t = 0.25$. (b) Dynamics for the  unstable wavenumber $k = 7$ at the instant $t = 0.3$. }\label{Fig9}
\end{figure}

\begin{figure}[h!]
\centering
{\includegraphics[scale=0.8]{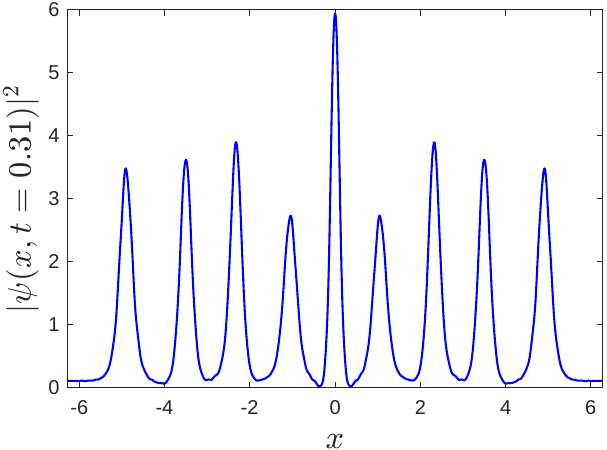}}
\caption{Dynamics for the system (\ref{EqA})-(\ref{EqB}) with $\epsilon = 10^{-4}$, parameter of dissipation $d_0 = 10^{-2}$, and unstable wavenumber $k = 3$. }\label{Fig10}
\end{figure}

The previous theorem demonstrates that, by selecting appropriate coefficients and initial conditions,  the general coupled NLS and the constant-coefficient coupled NLS, \emph{are equivalent at the initial stage of the dynamics}. In this scenario, we predict that the instability that is present in the latter will manifest itself in some manner in the former. For reliability of the hypothesis, and based on transformations (\ref{Transf1})-(\ref{Transf2}) and CW solutions (\ref{CWa})-(\ref{CWb}), we carried out numerical simulations for the system (\ref{EqA})-(\ref{EqB})  with initial conditions (perturbing such CW solutions)
\begin{equation}
  \psi(x,0) = \varphi(x,0) = 1 + \epsilon \cos(kx), \quad \quad 0 < \epsilon \ll 1.
\end{equation}
The coefficients used in the numerical experiments were $a(t) = e^{\cos t}$, $b(t) = \frac{1}{4} e^{-\cos t}\cos t $, $c(t) = 0,$ $d(t) = d_0 e^{2-\cos t}$, $g(t) = e^{2-\cos t},$ and $h(t) = -8 e.$ And then a solution of the Riccati system is
\begin{equation}
  \alpha(t) = -\frac{1}{4}e^{-\cos t}\sin t, \quad \quad \beta(t) = e^{1- \cos t}, \quad \quad \gamma(t) = \int_{0}^{t}e^{2-\cos s} ds, \quad \quad \mu(t) = e^{\cos t -1}. 
\end{equation}
The result established in Theorem 1 allows us to consider the instability region of the system without quadratic potential, that is, $|k| < \sqrt{32} e^{-d_0 \tau}$. Therefore, the goal is to investigate the existence of MI of the harmonic coupled NLS system for wavenumbers within or outside of this region. The results of the numerical simulations are reported in Figures \ref{Fig9}-\ref{Fig10}. Figure \ref{Fig9} depicts the system's response to energy supply ($d(t)<0,$ or equivalently $d_0 < 0$) and stable/unstable perturbations. In both situations, injection of energy into the system promotes the formation of a soliton train, as seen in snapshots. However, contrary to the case $b(t) = 0$, when the modulation instability produces solitons of almost equal amplitude, the quadratic potential generates solitons of different amplitudes (see also Figure \ref{Fig10}). In line with Theorem \ref{Theorem1}, the cause for instability in both circumstances is the expansion of the instability region, as the exponential's coefficient, $d_0$, is positive. 

Figure \ref{Fig10} displays the system's dynamics within a dissipative regime, that is, when $d(t)>0$. In this situation, modulational instability may emerge due to unstable wavenumbers and a small dissipation term. It is important to note that high dissipation rates tend to mitigate or control

 \noindent instability.
 
In conclusion, modulational instability can occur in the harmonic CNLS system and exhibit behavior similar to the non-harmonic case during the early stages of the dynamics, as long as suitable initial conditions are chosen, as established by Theorem \ref{Theorem1}.

\section{Conclusions and Final Remarks}
\label{Sect6}
This study delves into the intrinsic nonlinear dynamics of a newly proposed system of coupled nonlinear Schr\"odinger equations. Under specific conditions on the coefficients, the system is shown to support a rich set of solutions, including plane waves, Jacobi elliptic functions, and solitons. In particular, these solutions demonstrate robustness against small perturbations in the initial conditions, supporting their relevance for both physical modeling and numerical implementation. Similarly, the numerical simulations performed in this work indicate that specific external perturbations can lead the system into chaotic regimes. Therefore, implementing an appropriate control strategy is essential to suppress or regulate this unpredictable dynamic. Modulation instability was further analyzed, and the results suggest that proper adjustment of the dispersion and nonlinearity parameters can lead to the development of soliton trains within the system.

On the other hand, this work also contributes by developing analytical methods for deriving explicit solutions to nonlinear  Schr\"odinger-type models. These solutions, in turn, serve as valuable tools for assessing the stability and accuracy of the numerical techniques applied. Finally, it is hoped that the insights provided in this study will aid to understand a wide range of physical phenomena observed in real-world systems.

\begin{acknowledgement}
J.M. Escorcia thanks the Universidad EAFIT for the financial support provided for this research (Internal Project No. 12330022023).
\end{acknowledgement}

\section{Appendix: Solution of the Riccati system (\ref{Riccati1})-(\ref{Riccati4}) }
\label{Sect7}

In this appendix, we provide a solution of the Riccati system (\ref{Riccati1})-(\ref{Riccati4}), and we point out that all the formulas involved in such a solution have been verified previously in \cite{Ko-su-su}. A solution of the Riccati system with multiparameters is given by the following expressions \cite{CorderoSoto2008,Escorcia,Suslov12}:

\begin{equation}
\mu \left( t\right) =2\mu \left( 0\right) \mu _{0}\left( t\right) \left(
\alpha \left( 0\right) +\gamma _{0}\left( t\right) \right) ,  \label{mu}
\end{equation}%
\begin{equation}
\alpha \left( t\right) =\alpha _{0}\left( t\right) -\frac{\beta
_{0}^{2}\left( t\right) }{4\left( \alpha \left( 0\right) +\gamma _{0}\left(
t\right) \right) },  \label{alpha}
\end{equation}%
\begin{equation}
\beta \left( t\right) =-\frac{\beta \left( 0\right) \beta _{0}\left(
t\right) }{2\left( \alpha \left( 0\right) +\gamma _{0}\left( t\right)
\right) }=\frac{\beta \left( 0\right) \mu \left( 0\right) }{\mu \left(
t\right) }W\left( t\right) ,  \label{beta}
\end{equation}%
\begin{equation}
\gamma \left( t\right) = \gamma \left( 0\right) -\frac{\beta
^{2}\left( 0\right) }{4\left( \alpha \left( 0\right) +\gamma _{0}\left(
t\right) \right) },  \label{gamma}
\end{equation}%
subject to the initial arbitrary conditions $\mu \left( 0\right) ,$ $%
\alpha \left( 0\right) ,$ $\beta \left( 0\right) \neq 0,$ $\gamma (0)$. The functions $\alpha _{0}$, $\beta _{0}$, and
$\gamma _{0}$, are explicitly given by 
\begin{equation}
\alpha _{0}\left( t\right) =\frac{1}{4a\left( t\right) }\frac{\mu
_{0}^{\prime }\left( t\right) }{\mu _{0}\left( t\right) },  \label{alpha0}
\end{equation}%
\begin{equation}
\beta _{0}\left( t\right) =-\frac{W\left( t\right) }{\mu _{0}\left( t\right) 
},\quad W\left( t\right) =\exp \left( -\int_{0}^{t} c\left( s\right) \ ds\right) ,  \label{beta0}
\end{equation}%
\begin{equation}
\gamma _{0}\left( t\right) =\frac{1}{2\mu _{1}\left( 0\right) }\frac{\mu
_{1}\left( t\right) }{\mu _{0}\left( t\right) },  \label{gamma0}
\end{equation}%
Here $\mu _{0}$ and $\mu _{1}$
represent the fundamental solution of the characteristic equation subject to
the initial conditions $\mu _{0}(0)=0$, $\mu _{0}^{\prime }(0)=2a(0)\neq 0$
and $\mu _{1}(0)\neq 0$, $\mu _{1}^{\prime }(0)=0$.

\bibliography{references}
\bibliographystyle{unsrt} 


\end{document}